\documentclass[11pt]{amsart}
\usepackage{amsmath,amsfonts}
\usepackage{amssymb}
\usepackage{amscd}
\usepackage{amsthm}
\usepackage{yhmath}
\usepackage{subfigure}

\usepackage[all]{xy}
\usepackage{color}

\numberwithin{equation}{section}

\newtheorem{theorem}[equation]{Theorem}

\newtheorem{proposition}[equation]{Proposition}

\newtheorem{lemma}[equation]{Lemma}

\newtheorem{corollary}[equation]{Corollary}

\theoremstyle{remark}
\newtheorem{remark}[equation]{Remark}

\theoremstyle{definition}
\newtheorem{definition}[equation]{Definition}

\newtheorem{question}[equation]{Question}

\newtheorem{example}[equation]{Example}

\def\XXint#1#2#3{{\setbox0=\hbox{$#1{#2#3}{\int}$}
	\vcenter{\hbox{$#2#3$}}\kern-.5\wd0}}

\newcommand{\Ad}{\operatorname{Ad}}

\newcommand{\To}{\longrightarrow}

\newcommand{\N}{\mathbb N}

\newcommand{\C}{\mathbb C}

\newcommand{\R}{\mathbb R}

\newcommand{\G}{{\mathbb G}}
\newcommand{\g}{\mathfrak{g}}

\renewcommand{\L}{{\mathcal L}}

\newcommand{\diam}{\operatorname{diam}}

\newcommand{\id}{\operatorname{id}}

\newcommand{\vol}{\operatorname{vol}}
\renewcommand{\L}{{\mathcal L}}

\newcommand{\norm}[1]{\left\Vert#1\right\Vert}

\def\eps{\epsilon}

\begin{document}

\title[Smoothness  of subRiemannian  isometries]
{Smoothness  of subRiemannian  isometries}

\author{Luca Capogna}\address{Institute for Mathematics and its Applications, University of Minnesota, Minneapolis, MN 55455 \\ Department of Mathematical Sciences,
University of Arkansas, Fayetteville, AR 72701}\email{lcapogna@uark.edu}

\author{Enrico Le Donne}
\address{
Department of Mathematics and Statistics, University of Jyv\"askyl\"a, 40014 Jyv\"askyl\"a, Finland}
\email{enrico.ledonne@math.ethz.ch}

\keywords{Smoothness of isometries, Popp measure, subRiemannian geometry\\
LC is partially funded by NSF award  DMS 1101478}

\renewcommand{\subjclassname}{%
 \textup{2010} Mathematics Subject Classification}

\date{February 25, 2014}

\begin{abstract} 
  We    show that
  the group  of isometries 
   (i.e.,  distance-preserving homeomorphisms) 
  of an equiregular subRiemannian manifold is a finite-dimensional Lie group of smooth transformations.
  The proof is based on  a new PDE argument, in the spirit  of harmonic coordinates, establishing  that
  in an arbitrary subRiemannian manifold   
there exists an open dense subset   where
all isometries are smooth.
  
   \end{abstract}

\maketitle
\section{Introduction}

A classical result from 1939, due to Myers and Steenrod \cite{Myers-Steenrod}, states that any distance preserving homeomorphism (henceforth called isometry) between open subsets of smooth Riemannian manifolds is a smooth diffeomorphism. This regularity   opens the way for the use of analytical tools in the study of metric properties of manifolds. Since 1939 several different proofs for this result have been proposed for   Riemannian structures with varying degrees of smoothness. We refer the reader to the work of Palais \cite{Palais1}, Calabi and Hartman \cite{Calabi-Hartman}, and Taylor \cite{Taylor}. The latter   introduced a novel streamlined approach to the regularity issue, which (unlike many of the previous arguments) is not based on the study of the action of isometries on geodesics but rather relies on PDE. Namely, in \cite{Taylor}, Taylor notes that isometries preserve harmonic functions. Consequently,  the composition of the isometry with harmonic coordinates must be harmonic and as smooth as the ambient geometry.

Since the distance preserving property transfers easily to Gromov-Hausdorff limits for sequences of pairs of  Riemannian manifolds, in many applications it becomes relevant to study isometries between metric spaces obtained this way. In particular, in various settings including Mostow rigidity theory and the study of  analysis and geometry on boundaries of strictly pseudo-convex domains in $\C^n$,  the limit spaces fail to be Riemannian manifolds but have instead a natural subRiemannian structure. This gives rise to the problem of establishing regularity  of isometries in this broader setting. Since subRiemannian distance functions are generally not smooth with respect to the underlying differential structure, this regularity problem is quite different from its Riemannian counterpart.

Throughout this paper we define a subRiemannian manifold as a triplet  $(M,\Delta, g)$ where $M$ is a connected, smooth manifold of dimension $n\in \N$, $\Delta$ denotes a  subbundle of the tangent bundle $TM$ that bracket generates $TM$, and $g$ is a positive definite smooth, bilinear form defined on $\Delta$,
see \cite{Montgomery}.
 Analogously to the Riemannian setting, one can endow $(M,\Delta, g)$ with a metric space structure by defining the {\it Carnot-Caratheodory}  (CC) control distance: For any pair $x,y\in M$ set  \begin{multline} d(x,y)=\inf\{\delta>0 \text{ such that there exists a curve } \gamma\in C^\infty([0,1]; M) \text{ with endpoints }x,y \\\text{ such that }\dot\gamma\in \Delta (\gamma) \text{ and }|\dot\gamma|_g\le \delta\}.\notag\end{multline}
Curves whose velocity vector lies in $\Delta$ are called {\it horizontal}, their length  is defined in an obvious way. A horizontal curve is a  {\it geodesic} if it is  locally distance minimizing. A geodesic is {\it normal} if it satisfies a subRiemannian analogue of the geodesic equation. One of the striking features of subRiemannian geometry is that not all geodesics are normal (see \cite{Montgomery} and references therein). 

The first systematic study of  distance preserving homeomorphisms\footnote{Henceforth we will refer to distance preserving homeomorphisms as isometries and will specify if they are Riemannian or subRiemannian only when it is not completely clear from the context.} in the subRiemannian setting goes back to the groundbreaking   papers by Strichartz \cite{Strichartz, Strichartz2} and by Hamenst\"adt \cite{Hamenstadt}. In particular,  it follows from Hamenst\"adt's results that, assuming that all geodesics are normal,  isometries are smooth \cite[Theorem  6.2]{Hamenstadt}. If the underlying structure is a Carnot group (see Example \ref{carnot}), then  by virtue of the work of  Hamenst\"adt and of  Kishimoto \cite[Corollary 8.4]{Hamenstadt}, \cite[Theorem 4.2]{Kishimoto}, one has a stronger result:   Global isometries are compositions of group translations and group isomorphisms. More recently,   Ottazzi and the second named author  \cite[Theorem 1.1]{LeDonne-Ottazzi} have established a local version of this result, that applies to isometries between open subsets of Carnot groups and that does not rely on the Hamenst\"adt-Kishimoto's arguments.
Pansu's differentiability Theorem \cite{pansu} allows us to summarize
the previous results as follows.

\begin{theorem}[{\cite{LeDonne-Ottazzi}}]\label{H-theorem} Let $\G_1, \G_2$ be  Carnot groups and for $i=1,2$ consider  $\Omega_i\subset \G_i$  open sets. If $f:\Omega_1\to \Omega_2$ is an isometry then
$\G_1$ is isomorphic to $ \G_2$
and $f$ is the a  restriction to $\Omega_1$ of the composition of a group translation with   a Lie group automorphism.
\end{theorem} 

Smoothness  of isometries in this setting, where not all geodesics are normal, was settled for globally defined isometries by Kishimoto, using a theorem of Montgomery and Zippin  \cite[page 208, Theorem 2]{mz}, see also the detailed proof in \cite[Corollary  2.12]{LeDonne-Ottazzi},  while  in the local case (of mapping between open subsets of Carnot groups) it follows from the study of $1$-quasiconformal mappings by 
Cowling and the first named author \cite{Capogna-Cowling}.

The main results of the present paper represent an extension of the work mentioned above to the general subRiemannian setting, without the extra assumptions on the presence of  homogeneous structures or on the
absence of abnormal geodesics.

\begin{theorem}\label{generale}
Isometries between      subRiemannian manifolds are smooth 
in any open sets where the subRiemannian structure is equiregular.
In particular, there exists an  open dense  subset in which every   isometry is  smooth.

\end{theorem}
The proof of Theorem \ref{generale} proceeds along two steps of independent interest:  Theorem \ref{reg:theorem}
and Theorem \ref{Popp:preserved_intro}.

\begin{theorem}\label{reg:theorem}
Let $F:M\to N$ be an  isometry between two    subRiemannian manifolds.
If there exist  two $C^\infty$ volume forms
$\vol_M$ and $\vol_N$   such that
  $F_*\vol_M=\vol_N$,
then 
$F$ is a $C^\infty$ diffeomorphism.
\end{theorem}

A $C^\infty$ {\em volume form} on an $n$-dimensional manifold is the measure associated to a $C^\infty$ nonvanishing $n$-form on the manifold.
When $F:M\to N$ is a continuous map, the measure  $F_*\vol_M$ is defined by 
  $F_*\vol_M(A):=   \vol_M(F^{-1}(A))$, for all measurable sets $A\subseteq N$.
 
 The proof of   Theorem \ref{reg:theorem} is similar in spirit to   Taylor's approach to the regularity of Riemannian isometries via harmonic coordinates \cite{Taylor}. However, it appears  very difficult at present to establish existence of harmonic coordinates in the subRiemannian setting.
%
 To overcome this difficulty we observe that while isometries preserve harmonic functions (critical points of the Dirichlet energy \eqref{hor_energ}) they also preserve weak solutions of  nonhomogenous PDE in divergence form  built by polarization from the Dirichlet energy, see Corollary \ref{isom:preserve:sol}.
  Using this observation we show that the pull-back of {\it any} system of coordinates in the target manifold satisfies a subelliptic linear PDE for which   appropriate $L^p$ estimates  holds. Regularity  then follows through a bootstrap argument.

  Next we turn to the hypothesis of  Theorem \ref{reg:theorem}.
 The existence of smooth volume forms that are preserved by (metric) isometries is far from trivial. The natural choice for an isometric invariant measure on a manifold  is the top-dimensional  Hausdorff 
 measure $\mathcal{H}$
 (or spherical Hausdorff measure $\mathcal{S}$). In the Riemannian setting  $\mathcal{H}$ and $\mathcal{S}$ coincide with   the Riemannian volume measure and hence they readily provide  smooth
 measures that are isometric invariant. 
 A similar argument carries out to the Carnot group setting, where the Hausdorff measures are Haar measures and hence equal up to a constant to any left-invariant Riemannian measure.
However, for general 
{\em equiregular} subRiemannian manifolds   (see Definition \ref{equiregular})
the  counterpart of this statement is false. In \cite{Mitchell}, Mitchell shows that if $Q$ denotes the local homogeneous dimension (see  \eqref{dimension}), then the  density of Lebesgue measure with respect to the $Q$-dimensional   Hausdorff volume is bounded from above and below by positive constants. In their recent   paper \cite{Agrachev_Barilari_Boscain}, Agrachev, Barillari, and Boscain have studied the density function of any smooth volume form with respect to the spherical Hausdorff volume $\mathcal{S}$ in {\em equiregular} subRiemannian manifolds 
 and have proved that it is always continuous. Moreover,  if the topological dimension of $M$ is less or equal to $4$ this function is smooth. They also showed that for dimensions $5 $ and higher this is no longer true and exhibit examples where the density is $C^3$ but not $C^5$.  
 
 In view of this obstacle one has to find a different approach. Namely, we consider a canonical smooth volume form that is defined only in terms of the subRiemannian structure, the so called {\it Popp measure} (see \cite{Montgomery}) and show the following result.

\begin{theorem}\label{Popp:preserved_intro}
Let $F:M\to N$ be an  isometry between equiregular subRiemannian manifolds.
If $\vol_M$ and $\vol_N$ are the Popp measures on $M$ and $N$, respectively,
then 
 $F_*\vol_M=\vol_N$.
\end{theorem}

For smooth isometries a proof of the above result can be find  in \cite{Barilari-Rizzi}. The  novelty in 
Theorem~ \ref{Popp:preserved_intro} consists in the fact that we are considering  isometries with no smoothness assumptions. 

The proof  of Theorem \ref{Popp:preserved_intro} is based on the nilpotent tangent space approximation of the subRiemannian structure (see Section \ref{notation}) and on the detailed analysis of the induced map among such tangent spaces. Ultimately, the argument rests on Theorem \ref{H-theorem} and on Theorem \ref{abb}, which is a representation formula for Hausdorff measures recently established in \cite{Agrachev_Barilari_Boscain} and    \cite[Section 3.2, Remark 2]{ghezzi-jean}.

\begin{corollary}\label{cor:equiregular} Every $F:M\to N$   isometry between equiregular subRiemannian manifolds
is a $C^\infty$ diffeomorphism.\end{corollary}

Theorem \ref{generale} follows immediately from Corollary \ref{cor:equiregular} 
once we notice that, if $M$ is a subRiemannian manifold, then there exists an open dense  subset $U\subset M$  such that the subRiemannian structure restricted to  any connected component of $U$ is equiregular. 


We conclude this introduction 
with some consequences  of our regularity result, 
which in fact were our initial motivations.
Note that one can always extend an equiregular subRiemannian structure to a Riemannian one, i.e., there exists a Riemannian  tensor on $TM$ that coincides with $g$ when restricted to $\Delta$.
Riemannian extensions are also called Riemannian metrics  that {\it tame}  the  subRiemannian metric, \cite[Definition 1.9.1]{Montgomery}.
Hence, suppose that $F:M\to N$ is an isometry between two equiregular subRiemannian manifolds.
For any Riemannian extension $g_M$,  the push-forward  $g_N:=F_*g_M$  is well-defined  because of Corollary \ref{cor:equiregular} and tames 
 the subRiemannian structure on $N$. 
By construction, the map $F:(M,g_M)\to (N,g_N)$ is also a Riemannian isometry.

 
 It is less clear that one can find a single Riemannian extension for which a subRiemannian self-isometry is also a Riemannian isometry. This is the content of our main application: a structure theorem on the  group of subRiemannian self-isometries  ${\rm Isom}(M) $  and its correspondence to the group of Riemannian isometries associated to a  Riemannian metric that tames the subRiemannian one.

\begin{theorem}\label{app1}
If  $M$ is an equiregular subRiemannian manifold, then
\begin{enumerate}
\item[(i)]   ${\rm Isom}(M) $ admits a structure of  finite-dimensional Lie group;
\item[(ii)] for all compact subgroup  $K< {\rm Isom}(M)$, there exists a Riemannian extension $g_K$ on $M$ such that
  $K< {\rm Isom}(M,g_K)$.
\end{enumerate}
\end{theorem}
In particular, since the group of isometries that fixes a point $p\in M$,  denoted by $ {\rm Isom}_p(M)$, is compact (in view of Ascoli-Arzel\`a's theorem) 
then we have the following immediate consequence.

\begin{corollary}\label{E:extension}
If  $M$ is an equiregular subRiemannian manifold, then
for all $p\in M$, there exists a Riemannian extension $\hat g$ on $M$ such that
  ${\rm Isom}_p(M)< {\rm Isom}(M,\hat g)$.
\end{corollary} 


A further consequence is based on an argument in \cite{LeDonne-Ottazzi}. Just like Riemannian isometries, subRiemannian isometries are completely determined by the value and the horizontal differential at a point.
\begin{corollary}\label{horizontal:differential}
Let $M$ and $N$ be two connected equiregular subRiemannian manifolds.
Let $p\in M$ and let $\Delta$ be the horizontal bundle of $M$.
Let $f,g:M\to N$ be two isometries.
If $f(p)=g(p)$ and $df|_{\Delta_p}=dg|_{\Delta_p}$, then $f=g$.
\end{corollary}


\noindent{\it Acknowledgments}
The authors would like to thank E. Breuillard, U. Boscain, G. Citti, P. Koskela, U. Hamenst\"adt,   F. Jean, and A. Ottazzi for useful feedback and advice. Most of the work in this paper was developed while the authors were guests of the program {\it Interactions between Analysis and Geometry} at IPAM in the Spring 2013. The authors are very grateful to the program organizers M. Bonk, J. Garnett, U. Hamenst\"adt, P. Koskela and E. Saksman, as well to IPAM for their support. 
In particular, we thank Breuillard and Ottazzi 
for discussing with us  Corollary \ref{horizontal:differential}.


\section{Basic definitions and preliminary results}\label{notation}

%
%
%
%

Consider a subRiemannian manifold $(M,\Delta, g)$ and iteratively set $\Delta^1:=\Delta $,
and $\Delta^{i+1}:=\Delta^i+[\Delta^i, \Delta]$ for $i\in \N$. The bracket generating condition (also called {\it H\"ormander's finite rank hypothesis})
is expressed by the existence of $m\in \N$ such that,  for all $p\in M$, one has \begin{equation}\label{Hormander}\Delta_{p}^m= T_pM.\end{equation}

\begin{definition}\label{equiregular} A subRiemannian manifold $(M, \Delta, g)$ is {\it equiregular} if, for all $i\in \N$, the dimension of $\Delta^i_p$ is constant in $p\in M$. In this case, the {\it homogenous dimension}   is 
\begin{equation}\label{dimension} Q:=\sum_{i=1}^{m-1} i \,[\dim (\Delta_p^{i+1})-\dim (\Delta_p^i)].\end{equation}
\end{definition}

Consider the metric space $(M,d)$ where $(M,\Delta,g)$ is an equiregular subRiemannian manifold and $d$ is the corresponding control metric. 
As a consequence of Chow-Rashevsky
Theorem
 such a  distance is always finite  and induces on $M$ the original topology.
 As a result of Mitchell \cite{Mitchell}, the Hausdorff dimension of $(M,d)$ coincide with \eqref{dimension}.
 
 \begin{example}\label{carnot} {\it Carnot groups}  are important examples  of equiregular subRiemannian manifolds. 
 Classical references are 
 \cite{Folland75, Roth:Stein, Folland-Stein}.
 These are simply connected analytic  Lie groups $\G$ whose Lie algebras $\g$ has a stratification 
$\g=V^1\oplus ...\oplus V^m$ with$[V^1,V^j]=V^{1+j}$ and $[V^j,V^m]=0$ for all $j=1,...,m$. The corresponding subRiemannian structures are given by setting $\Delta$ to be the left-invariant subbundle associated to $V^1$
and by any choice of $g$ left-invariant metric on $\Delta$.
\end{example}
\subsection{Tangent cones and measures}\label{tangent}
The {\it tangent cone} of the metric space $(M,d)$ at a point $p\in M$ is the Gromov-Hausdorff limit 
$\mathcal{N}_p(M):=\lim_{t\to 0} (M,d/t,p)$.
 In view of Mitchell's work \cite{Mitchell} the metric space $\mathcal{N}_p(M)$ is described by
 the nilpotent approximation associated to the spaces  $\Delta^{i}_p$. In particular, $\mathcal{N}_p(M)$ is a Carnot group. 


%
%

For the definition of Popp measure on a subRiemannian manifold $(M, \Delta, g)$ we
refer the reader to \cite[Section 10.6]{Montgomery}. 
Here we only recall  that if $U\subset M$ is an open set upon which the subRiemannian structure is equiregular then 
 the Popp measure is a smooth volume form on $U$.
 

Following \cite{Agrachev_Barilari_Boscain}, we recall that,
given any $C^\infty$ volume form $\vol_M$ on $M$, e.g., the Popp measure,
we have a  $C^\infty$ volume form induced by $\vol_M$ on the nilpotent approximation  ${\mathcal N}_p(M)$, which
 we denote by  ${\mathcal N}_p(\vol_M)$.
Moreover, if $\vol_M$ is the Popp measure on $M$ and 
$\vol_{{\mathcal N}_p(M)}$ is the Popp measure on ${\mathcal N}_p(M)$,  then 
\begin{equation}\label{F2}
{\mathcal N}_p(\vol_M)   =   \vol_{{\mathcal N}_p(M)}.
\end{equation}
 
We shall use in a very crucial way the following representation formula for smooth volume measures. 
\begin{theorem}[{\cite[pages 358-359]{Agrachev_Barilari_Boscain},  \cite[Section 3.2]{ghezzi-jean}}]\label{abb}
Denote by $Q$ the Hausdorff dimension of an equiregular subRiemannian manifold $M$ and by ${\mathcal S}^Q_M$ the spherical Hausdorff 
measure on $M$.
Any $C^\infty$ volume form is related to ${\mathcal S}^Q$  by
\begin{equation}\label{F1}
\mathrm d \vol_M = 2^{-Q} {\mathcal N}_p(\vol_M) (B_{{\mathcal N}_p(M)} (e,1)) \mathrm d {\mathcal S}^Q_M.
\end{equation}
\end{theorem}
Here we denote by $B_{{\mathcal N}_p(M)} (e,1)$ the unit  ball in the metric space ${\mathcal N}_p(M)$ with center the identity element.

\begin{remark}\label{rmk:measure}
%

When read in charts, smooth volume forms are absolutely continuous with respect to the Lebesgue measure. 
The volume density is  smooth and bounded away from zero and infinity on compact sets.
Namely, for any smooth volume form $\vol$, one has $d\vol=\eta d {\mathcal L}^n$, with $\eta$ a smooth, strictly positive density function. Here and hereafter, by ${\mathcal L}^n$ we denote the Lebesgue measure.
 In particular, there is no ambiguity to say that some property holds almost everywhere without specifying with respect to which of the two possible measures: a volume form $\vol$ or Lebesgue measure.
\end{remark}

\subsection{Dirichlet energy}

Let  $X=\{X_1, \ldots, X_r \}$ be an orthonormal frame of $\Delta$ in a    subRiemannian manifold $(M,\Delta,g)$, and  $\vol$ a smooth volume measure on $M$.
Denote by ${\rm Lip}(M)$ the space of all Lipschitz functions  with respect
to the subRiemannian distance and recall the definition of the 
{\em horizontal gradient} 
\begin{equation}\label{hor_grad}
\nabla_{\rm H} u:= (X_1 u)X_1+\ldots +(X_r u)X_r, 
\end{equation}
defined a.e., as an $L^\infty_{\rm loc}(M)$ function,  for $u\in {\rm Lip}(M)$ (see \cite[Theorem 1.3]{Garofalo-Nhieu98}).
Given a compact set $K\subseteq M$,
we consider the {\em (horizontal Dirichlet) energy} (on $K$ and with respect to $\vol$) as
\begin{equation}\label{hor_energ}
E(u):=E_{{\rm H}, K}^{\vol}(u) :=
\int_K \norm{\nabla_{\rm H} u }^2\;\mathrm d\vol,
\end{equation}
and the associated quadratic form defined as, for all $u,v\in {\rm Lip}(M)$,
\begin{equation}\label{quad:form}
E(u,v)  :=E_{{\rm H},K} ^{\vol}(u,v):=
\int_K \langle \nabla_{\rm H} u ,  \nabla_{\rm H} v\rangle  \; \mathrm d\vol.
\end{equation}
Note that $2E(u,v)= E(u+v)- E(u) - E(v).$
The definitions above extend immediately to the larger class $W^{1,2}_{\rm H,loc}$, as defined in the next section. 



Regarding   the next result and an introduction to  the notion of (minimal) upper gradients, we refer the reader to   Hajlasz and Koskela's monograph.  
\begin{lemma}[{\cite[page 51 and Section 11.2]{HK}}]\label{minimal:gradient}
If $u\in {\rm  Lip}(M)$, then  $\norm{\nabla_{\rm H} u }$ coincides almost everywhere with the minimal upper gradient of $u$.
\end{lemma}

\begin{proposition}\label{isom:preserve:energy}
Let $F:M\to N$ be an isometry between two subRiemannian manifolds.
If $\vol_M$ and $\vol_N$ are two $C^\infty$ volume forms such that
 $F_*\vol_M=\vol_N$ then,
for all compact sets $K\subseteq M$,  
\begin{enumerate}
\item[(i)] 
$ E_{{\rm H},K}^{\vol_M}(u\circ F)=E_{{\rm H}, F(K)}^{\vol_N}(u)  , \qquad \forall u\in {\rm Lip}(N);$
\item[(ii)] $ E_{{\rm H},K} ^{\vol_M}(u\circ F, v\circ F)=E_{{\rm H}, F(K)}^{\vol_N}(u,v)  , \qquad \forall u,v\in {\rm Lip}(N).$
\end{enumerate}
\end{proposition}
\proof 
Since $F$ is distance  and volume preserving then it is an isomorphism between
the metric measure spaces $(M, d_M, \vol_M)$ and $(N,d_N,\vol_N)$. In particular, since  minimal upper gradients are a  purely metric measure notion, the map $F$ induces a one-to-one correspondence between minimal upper gradients between the two spaces.
Hence, by Lemma \ref{minimal:gradient}, being $\norm{\nabla_{\rm H} u } \circ F$ the pull back to $M$ of the  minimal upper gradient of $u $ in $(N,d_N,\vol_N)$, it coincides a.e. with  $\norm{\nabla_{\rm H} (u \circ F) }$, i.e.,  the minimal upper gradient of $u \circ F$ in $(M, d_M, \vol_M)$.
Consequently, $$\int_K \norm{\nabla_{\rm H} (u \circ F)}^2\;\mathrm d\vol_M 
=
\int_K \norm{\nabla_{\rm H} u }^2\circ F\;\mathrm d\vol_M = \int_{F(K)} \norm{\nabla_{\rm H} u}^2 \;\mathrm d\vol_N 
 .$$
 Since $E(\cdot)$ determines $E(\cdot, \cdot)$, the proof is concluded. 
\qed
%
%

We denote by ${\rm Lip}_{\rm loc}(M)$ the functions on $M$ that are Lipschitz on bounded sets.
We denote by 
${\rm Lip}_c(M)$ the Lipschitz functions on $M$ that have compact support.

\begin{corollary}\label{isom:preserve:sol}
If $u\in {\rm Lip}_{\rm loc}(N)$ and $g\in L^2(N)$ solve
 $$\int_N \langle \nabla_{\rm H} u ,  \nabla_{\rm H} v\rangle \; \mathrm d\vol_N= \int_N g v  \;\mathrm  d\vol_N  , \qquad \forall v\in {\rm Lip}_c(N),
$$
then the functions $\tilde u:= u\circ F$ and $\tilde g:= g\circ F$ solve
$$
\int_M \langle \nabla_{\rm H} \tilde u ,  \nabla_{\rm H} v\rangle\;  \mathrm d\vol_M= \int_M \tilde g v  \;  \mathrm d\vol_M,  \qquad \forall v\in {\rm Lip}_c(M).
$$
\end{corollary}
\proof
Note that, since $F$ preserves compact sets, it gives a one-to-one correspondence $v\mapsto v\circ F$ between ${\rm Lip}_c(N)$ and ${\rm Lip}_c(M)$.
Hence, any map in ${\rm Lip}_c(M)$ is of the form $v\circ F$.
Let $F(K)$ be a compact set containing the support of $v$.
Then 
\begin{eqnarray*}
\int_M \langle \nabla_{\rm H} \tilde u ,  \nabla_{\rm H} (v\circ F)\rangle\;   \mathrm d\vol_M &=&
 E_{{\rm H}, K} ^{\vol_M}(u\circ F, v\circ F)
 \\&=&
 E_{{\rm H}, F(K)}^{\vol_N}(u,v) 
 \\&=&  \int_N \langle \nabla_{\rm H} u ,  \nabla_{\rm H} v\rangle \; \mathrm d\vol_N
 \\&=& \int_N g\, v  \;  \mathrm d\vol_N
 \\
 &=&
 \int_M   (g\circ F )(v \circ F ) \; \mathrm d\vol_M.
\end{eqnarray*} 
\qed

\subsection{$L^p$ estimates for subLaplacians}

The regularity of isometries in Theorem \ref{reg:theorem} is based on a bootstrap argument, which eventually rests on  certain $L^p$ regularity estimates  established by Rothschild and Stein (see \cite[Theorem 18]{Roth:Stein}). In this section we recall these estimates in Theorem \ref{Rothschild-Stein}, along with the basic definitions of horizontal Sobolev spaces.

\begin{definition}\label{Sobolev spaces}
Let $X=\{X_1,...,X_r\}$ be a system of smooth vector fields in $\R^n$  satisfying H\"ormander's finite rank hypothesis \eqref{Hormander} and denote by $ {\mathcal L}^n$ the Lebesgue measure. For $k\in \N$ and for any multi-index $I=(i_1,...,i_k)\in \{1,...,r\}^k$ we define $|I|=k$ and $X^I u = X_{i_1}...X_{i_k} u$. For $p\in [1,\infty)$ we define the {\it horizontal Sobolev space} $W_{\rm H}^{k,p}(\R^n,  {\mathcal L}^n)$ associated to $X$  to be the space of all $u\in L^p(\R^n, {\mathcal L}^n)$  whose distributional derivatives $X^Iu$ are also in $L^p(\R^n,  {\mathcal L}^n)$ for all multi-indexes $|I|\le k$. This space can also be defined as the closure   of the space of $C^\infty_c(\R^n)$  functions   with compact support  with respect to the norm
\begin{equation}\label{sobolev}
\| u \|_{W_{\rm H}^{k,p}}^p := \| u\|_{L^p(\R^n, {\mathcal L}^n)}^p + \int_{\R^n} [\sum_{|I|=1}^k (X^I u)^2]^{p/2}  \;\mathrm d{\mathcal L}^n,
\end{equation}
see \cite{Garofalo-Nhieu98}, \cite{fssc-ms} and references therein.
A function $u\in L^p(\R^n,{\mathcal L}^n)$ is in the local Sobolev space $W^{k,p}_{\rm H,loc}(\R^n,{\mathcal L}^n)$ if, for any $\phi\in C^\infty_c(\R^n)$, one has $u\phi\in W_{\rm H}^{k,p}(\R^n,{\mathcal L}^n)$. 
\end{definition}
The definition above extends {\it verbatim} to a subRiemannian manifold $(M, \Delta, g)$ endowed  with a smooth volume measure $\vol_M$. The $W_{\rm H}^{k,p}$ Sobolev norm in this setting is defined through a   frame $X=\{X_1,...,X_r\}$ of $\Delta$ by
 \begin{equation}\label{sobolev-manifolds}\| u \|_{W_{\rm H}^{k,p}}^p := \| u\|_{L^p(M,\vol_M)}^p + \int_{M} [\sum_{|I|=1}^k (X^I u)^2]^{p/2}\, \mathrm d\vol_M.
\end{equation}

\begin{remark}\label{all the same}
 Although the Sobolev norm \eqref{sobolev-manifolds} 
depends on the specific frame and  smooth volume form chosen, 
the class  $ W^{k,p}_{\rm H,loc} (M, \vol)$ does not.
This is easily seen noticing that different choices of frames and smooth volume forms give rise to equivalent norms (on compact sets).
In particular, one can check whether a function is in $W^{k,p}_{\rm H,loc}$   through a smooth coordinate chart and using the Lebesgue measure, as in \eqref{sobolev}.
\end{remark}

\begin{theorem}[{\cite[Theorem 18]{Roth:Stein}}]\label{Rothschild-Stein}
Let $X_0, X_1,...,X_r$ be a system of smooth vector fields in $\R^n$  satisfying H\"ormander's finite rank hypothesis \eqref{Hormander}. Let 
\begin{equation}\label{sub-lap1}\mathcal L := \sum_{i=1}^r X_i^2 + X_0\end{equation} and consider a distributional solution to the equation $\mathcal Lu=f$ in $\R^n$. For every $k\in \N\cup \{0\}$ and $1<p<\infty$, if $f\in W_{\rm H}^{k,p}(\R^n,{\mathcal L}^n)$ then $u\in W^{k+2,p}_{\rm H,loc}(\R^n,{\mathcal L}^n)$.
\end{theorem} 

\begin{remark}\label{Lip is good}
If $u\in {\rm Lip}_{\rm loc}(\R^n)$ then $u\in W^{1,p}_{\rm H,loc}(\R^n,{\mathcal L}^n)$ for all $p\in [1,\infty)$. In fact, by virtue of 
\cite[Theorem 1.5]{Garofalo-Nhieu98} one has that the distributional derivatives $X_i u$ are in $L^\infty_{\rm loc}(\R^n)$.
\end{remark}
A standard argument, involving coordinate charts and cut-off functions, leads to a local version of the   Rothschild-Stein estimates  in any subRiemannian manifold.
For the reader's convenience we include the proof.
\begin{corollary}\label{regularity:sol}
Let $M$ be a subRiemannian manifold equipped with a $C^\infty$ volume form $\vol_M$ and let
$\ k\in \N$ and  $1\le p <\infty$.
If  $u\in W^{k,p}_{\rm H,loc}(M, \vol_M)\cap {\rm  Lip}_{\rm loc}(M)$  is
a solution of 
\begin{equation}\label{pde-intro}\int_M \langle \nabla_{\rm H} u ,  \nabla_{\rm H} v\rangle \; \mathrm d\vol_M= \int_M g v  \;  \mathrm d\vol_M   , \qquad \forall v\in {\rm Lip}_c(M).
\end{equation}
for some $g\in W^{k,p}_{\rm H,loc}(M,\vol_M)$, then $u\in W^{k+1,p}_{\rm H,loc}(M,\vol_M)$.
\end{corollary}


\proof
Consider an orthonormal frame $X=\{X_1,...,X_r\}$ of $\Delta$ and  a smooth coordinate chart $\phi: U\to V$, with $U\subset M$ and $V\subset \R^n$. Let 
$v \in  C^\infty_c(M)$ with  support in $U$. Without loss of generality we can assume that both $U$ and $V$ are relatively compact.
Setting
$$\tilde \Delta := \phi_* \Delta, \quad \tilde X:= \phi_* X, \quad 
\tilde u:= \phi_* u,\quad
\tilde v:= \phi_* v, \quad 
\tilde g:= \phi_* g, \quad
\tilde \vol_M:= \phi_* \vol_M,$$
then  through a  change of variable, \eqref{pde-intro}  reads as
$$\int_V  \sum_{i=1}^r \tilde X_i \tilde u  \tilde X_i \tilde v \,\mathrm d\tilde\vol_M= \int_V  \tilde g \tilde v  \,\mathrm d\tilde\vol_M.
$$
Since $\vol_M$ is a smooth volume measure, we can write 
$\mathrm d\tilde \vol_M= \eta  \mathrm d{\mathcal L}^n,$
for some $C^\infty$ function $\eta :V\to \R$ bounded away from zero.
Set  $$Y_1:=\sqrt{\eta} \tilde X_1, \ldots, Y_r:=\sqrt{\eta} \tilde X_r,$$ and note that the vector fields
$Y_1, \ldots,   Y_r$
are a  (orthogonal) frame for $\tilde \Delta$, and 
in particular they satisfy the H\"ormander finite rank condition \eqref{Hormander}. Since compactly supported smooth functions in $U$ are in a one to one correspondence with compactly supported smooth functions on $V$ under the action of $\phi$, it is immediate that
$$\int_V  \sum_{i=1}^r Y_i \tilde u  Y_i  v \,\mathrm d{\mathcal L}^n = \int_V \eta \tilde g v  \, \mathrm d{\mathcal L}^n
,\qquad \forall  v\in C^\infty_c(V),
$$
that is, $\tilde u$ is a distributional solution of the PDE
\begin{equation} \label{sub-lap2}\mathcal L \tilde u= \sum_{j=1}^r Y^*_j(Y_j(\tilde u))=\eta \tilde g.\end{equation}
Here, if $Y_i=\sum_{j=1}^n a_{ij}\partial_j$, then $Y_i^*$ denotes the {\it dual} vector (with respect to Lebesgue measure)
$Y_i^*=-\sum_{j=1}^n a_{ij} \partial_j -\sum_{j=1}^n \partial_j a_{ij}$. Consequently the left hand side of  \eqref{sub-lap2} is an operator $\mathcal L$ of the form \eqref{sub-lap1} with $X_0= \sum_{i=1}^r\sum_{j=1}^n \partial_j a_{ij}Y_i\in \tilde \Delta$. By virtue of Remark \ref{all the same} the functions $\tilde u,$ and $\tilde g\eta$ are in the Sobolev space $W^{k,p}_{\rm H,loc}(U, {\mathcal L}^n)$ corresponding to the H\"ormander vector fields $Y_1,...,Y_r$.  Next, in order to apply  Theorem \ref{Rothschild-Stein} we actually need to have the  right hand side in \eqref{sub-lap2} to belong to the smaller Sobolev space  $W_{\rm H}^{k,p}(\R^n, {\mathcal L}^n)$. To achieve this we consider for every $p\in U$ a cut-off function $\varphi\in C^{\infty}_c(U)$ which is identically one near $p$. It will  be then enough to prove that $u\varphi \in W^{k+1,p}_{\rm H}(\R^n,{\mathcal L}^n)$. Evaluating $\mathcal L (u\varphi)$ yields
$$\mathcal L (u\varphi)= \varphi \mathcal L u + \sum_{i=1}^r b_i Y_i u + b u,$$
where $b_i,b\in C^{\infty}_c(U)$. Since $b \tilde u, \varphi\mathcal L \tilde u\in W_{\rm H}^{k,p}(\R^n, {\mathcal L}^n)\subset W_{\rm H}^{k-1,p}(\R^n, {\mathcal L}^n)$ and $b_i Y_i u\in W_{\rm H}^{k-1,p}(\R^n, {\mathcal L}^n)$ then the right hand side of the latter equation is in
$ W_{\rm H}^{k-1,p}(\R^n, {\mathcal L}^n)$  for all $p\in [1,\infty)$. The conclusion now follows from Theorem \ref{Rothschild-Stein} and from applying again Remark \ref{all the same}.
\qed


Smoothness of Sobolev functions is guaranteed by a Morrey-Campanato type embedding in an appropriate H\"older class.
\begin{definition}\label{Hoelder}
Let $X=\{X_1,...,X_r\}$ be a system of smooth vector fields in $\R^n$  satisfying H\"ormander's finite rank hypothesis \eqref{Hormander} inducing a control distance $d$. 
For $\alpha\in (0,1)$ define the (Folland-Stein) H\"older class 
$$C^\alpha_{\rm loc}(\R^n):=\left.\big\{u:\R^n\to \R  \;\right|\; \forall \text{ compact   } K , \exists C>0 :  \forall x,y\in K, 
{u(x)-u(y)}\le C {d^\alpha(x,y)}\big\}.$$
\end{definition}

For any system of smooth vector fields satisfying H\"ormander's finite rank hypothesis and for any
bounded open set $\Omega\subset \R^n$, one can find  constants  $Q=Q(\Omega,X)>0$, the so-called {\it local homogenous dimension} relative to $X$ and $\Omega$,  $C=C(\Omega,X)>0$ and $R=R(\Omega,X)>0$, such that, for all metric balls centered in $\Omega$ with  radius  less than $R$,
one has a lower bound on the Lebesgue volume of these balls in terms of a power of their radius, i.e., ${\mathcal L}^n(B)\ge C \diam(B)^Q$. The quantity $Q$ arises out of the celebrated doubling property established by  Nagel, Stein and Wainger  in \cite{nagelstwe} and plays an important role in the theory of Sobolev spaces in the subRiemannian setting. In  case $X$ is an equiregular distribution then  $Q$ does not depend on $\Omega$ and coincides with the homogenous dimension defined earlier.
 By virtue of Jerison's Poicar\'e inequality \cite{Jerison} one has that, if $u\in W^{1,p}_{\rm H, loc}(\R^n)$ and  $p>Q>1$, then, for any $d$-metric ball $B\subset \R^n$, one has 
 \begin{equation}\label{Campanato}
 \frac{1}{{\mathcal L}^n(B)}\int_B |u-u_B|^Q d{\mathcal L}^n\le C\text{diam} (B)^{Q(\frac{p-Q}{p})},
 \end{equation}
where $u_B$ denotes the average of $u$ over $B$ and $C>0$ is a constant depending only on $p,Q, X,$ and $\norm{u}_{L^p(B)}$. 
Then a simple argument, see for instance \cite{Lu98} and \cite[Theorem 6.9]{Capogna}, 
starting from \eqref{Campanato} leads to 
  the following Morrey-Campanato type embedding.
\begin{proposition}\label{Morrey}
Let $X=\{X_1,...,X_r\}$ be a system of smooth vector fields in $\R^n$  satisfying H\"ormander's finite rank hypothesis \eqref{Hormander}.  
For every bounded open set $\Omega\subset \R^n$ there exists a positive number $Q$ depending only on $n,X$ and $\Omega$ such that,  if $p>Q$, then $$W^{1,p}_{\rm H,loc}(\Omega,{\mathcal L}^n)\subset C^\alpha_{\rm loc}(\Omega)$$ with $\alpha=(p-Q)/p$.
\end{proposition}

%
%

\section{Proofs}
\subsection{Proof of Theorem \ref{Popp:preserved_intro}}



 


We begin by establishing the result in the special case of Carnot groups.
Assume that $M,N$ are two Carnot groups and let  $G: M \to N$ be an isometry fixing the identity.
In view of Theorem \ref{H-theorem} the map $G$ is a group isomorphism.
Since on Carnot groups the Popp measures are left-invariant,
$G_* (\vol_M )$ is a Haar measure, hence it is a constant multiple of 
$\vol_N.$
 Since $G$ is an isometry, the constant is one.

Let us go back to the general case of equiregular subRiemannian manifolds.
Recall the definition of nilpotent approximations ${\mathcal N}_p(M) $ from Section \ref{tangent}.
We observe that, for {\em every} $p\in M$, 
since $(M,p)$ and $(N,F(p))$ are isometric as pointed metric spaces, then, by uniqueness of Gromov-Hausdorff limits, 
there exists a (not necessarily unique) isometry 
$G_p :  {\mathcal N}_p(M)\to {\mathcal N}_{F(p)}(N)$, fixing the identity. 
In particular, 
the unit balls of the tangents are the same, i.e., 
\begin{equation}\label{F3}
G_p (B_{{\mathcal N}_p(M)} (e,1))
 =
 B_{{\mathcal N}_{F(p)}(N)} (e,1).
\end{equation}
In view of the argument at the beginning of the proof, we have that 
the Popp measures are the same, i.e., 
\begin{equation}\label{F4}
(G_p)_* (\vol_{{\mathcal N}_p(M)} ) = 
\vol_{{\mathcal N}_{F(p)}(N)}.
\end{equation}

In addition, the spherical Hausdorff measures are preserved by isometries. Namely, we have
\begin{equation}\label{F5}
F_*{\mathcal S}^Q_M =
{\mathcal S}^Q_N.
\end{equation}

Therefore, for all measurable sets $A\subseteq N$, we have
\begin{eqnarray*}
F_* \vol_M(A) & = & \vol_M (F^{-1}(A))\\
(\text{from \eqref{F1})\quad\qquad} &=& \dfrac{1}{2^Q} \int_{F^{-1}(A)}    {\mathcal N}_p(\vol_M) (B_{{\mathcal N}_p(M)} (e,1)) \,\mathrm d {\mathcal S}^Q_M(p)\\
(\text{from \eqref{F3})\quad\qquad} &=& \dfrac{1}{2^Q} \int_{F^{-1}(A)}    {\mathcal N}_p(\vol_M) \left(
G_p ^{-1} (B_{{\mathcal N}_{F(p)}(N)} (e,1) )\right)
\mathrm d {\mathcal S}^Q_M(p)\\
 &=& \dfrac{1}{2^Q} \int_{F^{-1}(A)}     (G_p)_* \,{\mathcal N}_p(\vol_M) \left(
 B_{{\mathcal N}_{F(p)}(N)} (e,1) \right)
\mathrm d {\mathcal S}^Q_M(p)\\
 (\text{from \eqref{F2})\quad\qquad} &=&
  \dfrac{1}{2^Q} \int_{F^{-1}(A)}     (G_p)_* \vol_{{\mathcal N}_p(M)} \left(
 B_{{\mathcal N}_{F(p)}(N)} (e,1) \right)
\mathrm d {\mathcal S}^Q_M(p)\\
 (\text{from \eqref{F4})\quad\qquad} &=&  
 \dfrac{1}{2^Q}   \int_{F^{-1}(A)}     \vol_{{\mathcal N}_{F(p)}(N)} \left(
 B_{{\mathcal N}_{F(p)}(N)} (e,1) \right)
\mathrm d {\mathcal S}^Q_M(p)\\
  (\text{$q=F(p)$ and from \eqref{F5})\quad\qquad} 
&=&
 \dfrac{1}{2^Q} \int_{A}  \vol_{ {\mathcal N}_q(N) }(B_{{\mathcal N}_q(N)} (e,1))\, \mathrm d {\mathcal S}^Q_N(q)\\
 (\text{from \eqref{F2})\quad\qquad} &=&
 \dfrac{1}{2^Q} \int_{A}   {\mathcal N}_q(\vol_N) (B_{{\mathcal N}_q(N)} (e,1)) \,\mathrm d {\mathcal S}^Q_N(q)\\
(\text{from \eqref{F1})\quad\qquad}  &=& \vol_N(A).
\end{eqnarray*}
The proof is concluded.
\qed

\begin{remark}
We observe that, in view of Margulis-Mostow's version of Rademacher Theorem
 \cite{Margulis-Mostow},
in the proof above there is a natural choice for  the maps $G_p$ for almost every $p$.
\end{remark}

\subsection{The iteration step}
Theorem \ref{reg:theorem} will be proved bootstrapping the following two results.
\begin{proposition}\label{Smoothnes_composition}
Let $M$ and $N$ be two subRiemannian manifolds equipped with $C^\infty$ volume measures $\vol_M$ and $\vol_N$.
Let $F:M\to N$ be an isometry such that  $F_*\vol_M=\vol_N$.
Let $g\in C^\infty(N)$.
Assume that
\begin{enumerate}
\item[(i)]  $g\circ F\in W^{k,p}_{\rm H,loc}(M,\vol_M)$, for  some $k\in \N$ and some $1\le p <\infty$,  and
\item[(ii)] The function $u\in  W^{k,p}_{\rm H,loc}(M,\vol_M)\cap {\rm  Lip}(N,\R)$ is a solution of $$\int_N \langle \nabla_{\rm H} u ,  \nabla_{\rm H} v\rangle \;  \mathrm d\vol_N= \int_N g \,v  \;\mathrm  d\vol_N , \qquad \forall v\in {\rm Lip}_c(N).
$$
\end{enumerate}
Then $u\circ F\in W^{k+1,p}_{\rm H,loc}(M, \vol_M)$ and solves
$$
\int_M \langle \nabla_{\rm H}  (u \circ F) ,  \nabla_{\rm H} v\rangle\;  \mathrm d\vol_M= \int_M  (g\circ F) \,v  \; \mathrm d\vol_M,  \qquad \forall v\in {\rm Lip}_c(M).
$$
\end{proposition}
\proof The last statement is  an immediate consequence of Corollary \ref{isom:preserve:sol}. Since $g\circ F\in W^{k,p}_{\rm H,loc}(M,\vol_M)$ then the conclusion follows from 
Corollary \ref{regularity:sol}.
\qed

\begin{proposition}\label{gF} 
Let $U$ and $V$ be two open sets of $\R^n$ and $X=\{X_1,...,X_r\}$ a frame of smooth H\"ormander vector fields in $U$. For $k\in \N\cup \{0\}$ and $1\le p <\infty$, denote by $ W^{k,p}_{\rm H,loc}(U,{\mathcal L}^n)$ the corresponding Sobolev space.
Let $F=(F_1,...,F_n) :U\to V$ be a map such that,  for all  $i=1,\ldots, n$, one has  $F_i \in W^{k,p}_{\rm H,loc}(U,{\mathcal L}^n)$. If  $g\in C^\infty (V)$, then
 $g \circ F \in W^{k,p}_{\rm H,loc}(U,{\mathcal L}^n)$.
\end{proposition}
\proof
For all $i=1,...,n$ and $j=1,...,r$, we have $X_j F_i \in W^{k-1,p}_{\rm H,loc}(U,{\mathcal L}^n)$. Using the scalar chain rule for Sobolev functions one obtains 
\begin{equation}\label{rhs} X_j(g\circ F)=  X_j(g(F_1,\ldots, F_n))=
\partial_1 g X_j F_1+\ldots \partial_n g X_j F_n,\end{equation} a.e. in $U$.
Since $\partial_i g\in C^\infty $, then the right hand side of \eqref{rhs} is also in $W^{k-1,p}_{\rm H,loc}(U,{\mathcal L}^n)$, concluding the proof.
\qed

\subsection{Proof of  Theorem \ref{reg:theorem}}
We work in (smooth) coordinates, so that 
the isometry can be considered as a map   $F :U\to V$  with $U$ and $V$  two open sets of $\R^n$.
Denote by  $x_1, \ldots, x_n$ the  coordinate functions in $V$.
The volume form of the second manifold can be written as 
$$\mathrm d\vol_N=\eta \mathrm d{\mathcal L}^n,$$
for some
$C^\infty$ function $\eta :V\to \R$ bounded away from zero on compact sets.

For each $i=1,\ldots, n$, set
  $g_i:= \eta^{-1} \nabla_{\rm H}^*(\eta  \nabla_{\rm H} x_i)$. We have
  $g_i\in C^\infty(V)$ and
 $$\int_V \langle \nabla_{\rm H} x_i ,  \nabla_{\rm H} v\rangle \eta \; \mathrm d{\mathcal L}^n
 = \int_V g_i v  \eta\;\mathrm  d{\mathcal L}^n  , \qquad \forall v\in {\rm Lip}_c(V).
$$
Notice that, for all $i=1,...,n$, one has that $x_i\circ F$ and $g_i\circ F$ are  Lipschitz function with respect to the control distance.
 Recalling Remark \ref{Lip is good}, we observe that as a consequence
$x_i\circ F, g_i\circ F \in W^{1,p}_{\rm H,loc}(U)$, for all $1\le p < \infty$.
Using Proposition \ref{Smoothnes_composition} with $k=1$, we have that 
$F_i:=x_i\circ F \in  W^{2,p}_{\rm H,loc}(U,{\mathcal L}^n) $.
By Proposition \ref{gF}, it follows that 
$g \circ F  \in W^{2,p}_{\rm H,loc}(U,{\mathcal L}^n)  $.
Iterating Proposition \ref{Smoothnes_composition}  and  Proposition \ref{gF}, one obtains $F_i \in W^{k,p}_{\rm H,loc}(U,{\mathcal L}^n)  $ for all $k\in \N$. Invoking 
Proposition \ref{Morrey} we finally conclude 
$ F_i  \in C^\infty(U) $, for all $i=1,...,n$.
\qed

\subsection{Proof of Theorem \ref{app1}}
%
%
Since ${\rm Isom}(M)$ is locally compact  by Ascoli-Arzel\`a's theorem, from a result of Montgomery and Zippin \cite[page 208, Theorem 2]{mz}, we have (i).
To show (ii), we fix an auxiliary Riemannian extension $\tilde g$ on $M$ and define
$$g=\int_H F_* \tilde g  \; \mathrm d \mu_H(F),$$
where $\mu_H$ is a probability Haar measure on $H$. 
Notice that $\mathcal G:=\{ F_*g : F\in H\}$ is a compact set of Riemannian tensors extending the subRiemannian structure on $M$.
In particular, in local coordinates,  all $g'\in \mathcal G$ can be represented with a matrix with uniformly bounded eigenvalues. Hence, $g$ is a Riemannian tensor, which is $H$-invariant, by linearity of integrals. 
\qed

\subsection{Proof  of Corollary \ref{horizontal:differential}}
Consider $F:=g^{-1}\circ f$. So $F$ is an isometry of $M$ that fixes $p$ and such that $dF|_{\Delta_p}=\id_{\Delta_p}$.
By Corollary \ref{E:extension}, there exists a Riemannian metric $\hat g$ for which $F$ is a Riemannian isometry.
Since $\Delta$ is bracket generating, one shows (see   \cite{LeDonne-Ottazzi}) that $dF|_{T_pM}=\id_{T_pM}$ and hence $F=\id_{M}$.
\qed


 \bibliography{general_bibliography}
\bibliographystyle{amsalpha}

\end{document}